\documentclass[12pt]{article}
\usepackage{amsmath}
\usepackage{amsfonts}
\usepackage{mathrsfs}
\usepackage{amssymb}
\usepackage{indentfirst}
\usepackage[usenames]{color}
\usepackage[numbers,sort&compress]{natbib}
\usepackage{subfigure}

\def\cl{\centerline}
\def\ni{\noindent}

\def\BB{{\cal B}(q)}
\def\sp{{\rm{span}}}
\def\deg{{\rm deg}}

\def\Z{\mathbb{Z}}

\def\C{\mathbb{C}}

\def\QED{\hfill$\Box$}

\def\Vir{{\rm {Vir}}}

\def\Uh0{\mathcal{U}(\mathfrak h_0)}
\def\Bp0{\mathcal{B}(p,0)}
\def\bp1{\mathcal{B}(p,1)}

\numberwithin{equation}{section}
\newtheorem{theo}{Theorem}[section]

\newtheorem{lemm}[theo]{Lemma}

\newtheorem{prop}[theo]{Proposition}
\newtheorem{clai}{Claim}

\newtheorem{exam}[theo]{Example}

\begin{document}
\begin{center}
{\bf\Large Representations of non-finitely graded Lie algebras related to Virasoro algebra}
\footnote{
$^{*}$Corresponding author: chgxia@cumt.edu.cn (C.~Xia).}
\end{center}
\vspace{6pt}

\cl{Chunguang Xia$^{\,*,\,\dag,\,\ddag}$, Tianyu Ma$^{\,\dag}$, Xiao Dong$^{\,\dag}$, Mingjing Zhang$^{\,\dag}$}
\vspace{6pt}

\cl{\small $^{\,\dag}$School of Mathematics, China University of Mining and Technology, Xuzhou 221116, China}
\vspace{6pt}

\cl{\small $^{\,\ddag}$School of Mathematical Sciences, Xiamen University, Xiamen 361005, China}
\vspace{6pt}

\cl{\small E-mails: chgxia@cumt.edu.cn, tyma@cumt.edu.cn, xdong@cumt.edu.cn, zhangmj@nju.edu.cn}
\vspace{6pt}

\footnotesize
\noindent{{\bf Abstract.}
In this paper, we study representations of non-finitely graded Lie algebras $\mathcal{W}(\epsilon)$ related to Virasoro algebra, where $\epsilon = \pm 1$.
Precisely speaking, we completely classify the free $\mathcal{U}(\mathfrak h)$-modules of rank one over $\mathcal{W}(\epsilon)$,
and find that these module structures are rather different from those of other graded Lie algebras.
We also determine the simplicity and isomorphism classes of these modules.

\ni{\bf Key words:} Virasoro algebra; non-finitely graded Lie algebras; free $\mathcal{U}(\mathfrak h)$-module; simplicity

\ni{\it Mathematics Subject Classification (2020):} 17B10; 17B65; 17B68}

\small
\section{\large{Introduction}}

It is well-known that the Virasoro algebra, denoted $\widehat{{\rm Vir}}$, is an infinite dimensional Lie algebra with basis $\{L_i,\,c\,|\,i\in\Z\}$ satisfying
$$
[L_i, L_j]=(j-i)L_{i+j}+\delta_{i+j,0}\frac{i^3-i}{12}c,\quad [L_i, c]=0.
$$
Here $c$ is the center of $\widehat{{\rm Vir}}$. We denote by ${\rm Vir}$ the centerless Virasoro algebra.
Clearly, ${\rm Vir}$ is a $\Z$-graded Lie algebra with Cartan subalgebra $\mathfrak{h}_{\rm Vir}=\C L_0$.
The representation theories of $\widehat{{\rm Vir}}$ and ${\rm Vir}$ have been extensively and deeply studied due to its importance in mathematics and physics,
see, e.g., the survey in \cite{IK2010}.
One of the most important work in weight module theory is the well-known Mathieu's theorem \cite{M1992} on classification of Harish-Chandra modules,
which was conjectured earlier by Kac \cite{K1981}.

Recently, an interesting non-weight module problem for a Lie algebra $\mathfrak g$, defined by an ``opposite condition'' relative to weight modules,
was proposed by Nilsson \cite{N2015,N2016}.
Let $\mathfrak h$ be the Cartan subalgebra of $\mathfrak g$. Denote by $\mathcal{U}(\mathfrak h)$ the universal enveloping algebra of $\mathfrak h$.
This kind of non-weight modules, referred to as {\it free $\mathcal{U}(\mathfrak h)$-modules} since the action of $\mathfrak h$ is required to be free,
was constructed first for $\mathfrak{sl}_n$ by Nilsson \cite{N2015} and independently by Tan and Zhao \cite{TZ2018}.
In another paper, Tan and Zhao \cite{TZ2015} proved that any free $\mathcal{U}(\mathfrak{h}_{\rm Vir})$-module of rank one over ${\rm Vir}$
is isomorphic to $\Omega_\Vir(\lambda,\alpha)$ (cf.~\eqref{def-vir-mod}) for some $\lambda\in\C^*$ and $\alpha\in\C$.
Following \cite{N2015,N2016,TZ2015,TZ2018}, free $\mathcal{U}(\mathfrak h)$-modules, especially for finitely graded Lie algebras containing Virasoro subalgebra,
have been extensively studied in recent years, see, e.g., \cite{CC2015,CG2017,HCS2017,CH2019,HCL2020,W2021} and the references therein.
It has been realized that interesting free $\mathcal{U}(\mathfrak h)$-module results will appear under the assumption that
$\mathfrak h$ is not a Cartan subalgebra (e.g.,~\cite{CG2017,HCS2017,N2023}). The latest work in this direction is Nilsson's work \cite{N2023},
where he chose $\mathfrak h$ to be the nilradical of a maximal parabolic subalgebra of $\mathfrak{sl}_n$.

Let $\epsilon=\pm 1$. In this paper, we focus on a class of non-finitely graded Lie algebras $\mathcal{W}(\epsilon)$
with basis $\{L_{i,m}\, |\,i\in\Z, m\in\Z_+\}$ and relations
\begin{equation}\label{brackets-for-w}
[L_{i,m},L_{j,n}]=(j-i)L_{i+j,m+n}+\epsilon(m-n)L_{i+j,m+n-\epsilon}.
\end{equation}
Note that the subspace spanned by $\{L_{i,0}\in \mathcal{W}(\epsilon)\,|\, i\in\Z\}$ is isomorphic to the centerless Virasoro algebra ${\rm Vir}$.
For simplicity, we refer to $\mathcal{W}(\epsilon)$ as {\it non-finitely graded Virasoro algebras}.
The case for $\epsilon=1$ was first constructed in \cite{SXZ2000,X2000} as the simple Lie algebra $\mathit{W}(0,1,0;\Z)$ of Witt type.
The case for $\epsilon=-1$ was naturally realized in \cite{CHS2014} as the non-simple Lie algebra ${W}(\Z)$ consisting of smooth function
$L_{i,m}:=-(1+t)^{-m} e^{-it}\in C^{\infty}_{[0,+\infty)}$ under bracket $[f,g]=fg'-f'g$ for $f,g\in{W}(\Z)$.

For non-finitely graded Lie algebras with Virasoro subalgebra, motivated by Mathieu's work \cite{M1992}, it is natural to consider the classification of the so-called quasifinite modules,
such as the $W$-infinity algebra $\mathcal{W}_{1+\infty}$ \cite{KR1993}, Lie algebras of Weyl type \cite{S2003}, Lie algebras of Block type \cite{SXX2012,SXX2013} and so on.
However, for the Lie algebra $\mathcal{W}(\epsilon)$, although it admits a natural principal $\Z$-gradation
$\mathcal{W}(\epsilon)=\oplus_{i\in\Z}\mathcal{W}(\epsilon)_i$ with $\mathcal{W}(\epsilon)_i=\sp\{L_{i,m}\,|\,m\in\Z_+\}$,
the zero-graded part $\mathcal{W}(\epsilon)_0$ is not a commutative subalgebra and worse it does not contain Cartan subalgebra.
This leads to that one cannot cope with the quasifinite modules over $\mathcal{W}(\epsilon)$.
Fortunately, in this paper, we find that one can develop the representation theory of $\mathcal{W}(\epsilon)$ by studying free $\mathcal{U}(\mathfrak h)$-modules
with $\mathfrak h=\C L_{0,0}$ (although it is not a Cartan subalgebra).
To best of our knowledge, up to now, there are few work on free $\mathcal{U}(\mathfrak h)$-modules over non-finitely graded Lie algebras,
see \cite{CY2018,GWL2021,WXZ2023} for some attempt on loop Virasoro algebra and some Lie algebras of Block type.

We surprisingly find that the module structures of $\mathcal{W}(\epsilon)$ are much more complicated than that of the Virasoro algebra \cite{TZ2015},
and rather different from those of other non-finitely graded Lie algebras \cite{CY2018,GWL2021,WXZ2023}.
Our first main result for $\mathcal{W}(1)$ is Theorem~\ref{thm-1}, in the proof of which the following combinatoric formula: for $m<n$,
\begin{equation}\label{formula}
\binom{n-1}{m} + \binom{n-1}{m-1} = \binom{n}{m}
\end{equation}
will be used frequently.
Our second main result for $\mathcal{W}(-1)$ is Theorem~\ref{thm-2}, for which we need more technical analysis.
We also would like to point out that our techniques used here may be applied to
analogous problems of non-finitely graded Lie algebras which are closely related to $\mathcal{W}(\epsilon)$.
This is also our motivation for writing this paper.

This paper is organized as follows. In Section 2, we recall the classification of free $\mathcal{U}(\mathfrak{h}_{\rm Vir})$-modules of rank one over ${\rm Vir}$,
and present an elementary result on a number sequence.
In Sections 3 and 4, we classify the free $\mathcal{U}(\mathfrak h)$-modules of rank one over $\mathcal{W}(1)$ and $\mathcal{W}(-1)$, respectively.
Along the way, we also determine the simplicity and isomorphism classes of these modules.

\section{\large{Preliminaries}}

Throughout this paper, we use $\Z$, $\Z_+$, $\C$ and $\C^*$ to
denote the sets of integers, nonnegative integers, complex numbers and nonzero complex numbers, respectively.
We work over the complex field $\C$.
In this section, we recall the classification of the free $\mathcal{U}(\mathfrak{h}_{\rm Vir})$-modules of rank one over the centerless Virasoro algebra ${\rm Vir}$.
We also present an elementary result on a number sequence, which will be used in Section~4.

\subsection{Free $\mathcal{U}(\mathfrak{h}_{\rm Vir})$-modules of rank one over ${\rm Vir}$}

Recall that $\mathfrak{h}_{\rm Vir}=\C L_{0}$ is the Cartan subalgebra of ${\rm Vir}$.
Let $\lambda\in\C^*$ and $\alpha\in\C$. If we define the action of ${\rm Vir}$ on the vector space of polynomials in one variable
$$
\Omega_{{\rm Vir}}(\lambda,\alpha):=\C[t]
$$
by
\begin{equation}\label{def-vir-mod}
L_i\cdot f(t)=\lambda^i(t-i\alpha)f(t-i),
\end{equation}
where $i\in\Z, f(t)\in\C[t]$,
then $\Omega_{{\rm Vir}}(\lambda,\alpha)$ becomes a ${\rm Vir}$-module.
Equivalently, we can rewrite \eqref{def-vir-mod} in a more simple form by restricting the action of $L_i$ on monomials:
 \begin{equation}\label{def-vir-mod-simple}
L_i\cdot t^k=\lambda^i(t-i\alpha)(t-i)^k.
\end{equation}
These modules firstly appeared in \cite{LZ2014} as quotient modules of fraction Virasoro modules.
From \eqref{def-vir-mod-simple}, we see that $\Omega_{{\rm Vir}}(\lambda,\alpha)$ is free of rank one when restricted to $\mathfrak{h}_{\rm Vir}$.
In fact, we have the following classification result \cite{TZ2015}.

\begin{lemm}\label{lemm-12}
Any free $\mathcal{U}(\mathfrak{h}_{\rm Vir})$-module of rank one over ${\rm Vir}$ is isomorphic to $\Omega_{{\rm Vir}}(\lambda,\alpha)$
defined by \eqref{def-vir-mod-simple} for some $\lambda\in\C^*$ and $\alpha\in\C$.
\end{lemm}

The simplicity and isomorphism classification of $\Omega_{{\rm Vir}}(\lambda,\alpha)$  are as follows \cite{TZ2015}.

\begin{lemm}\label{lemm-13}\parskip-3pt
\begin{itemize}\parskip-3pt
\item[{\rm(1)}] $\Omega_{\rm Vir}(\lambda,\alpha)$ is simple if and only if $\alpha\neq0$.
\item[{\rm(2)}] $\Omega_{\rm Vir}(\lambda,0)$ has a unique proper submodule
$t\Omega_{\rm Vir}(\lambda,0)\cong\Omega_{\rm Vir}(\lambda,1)$, and $\Omega_{\rm Vir}(\lambda,0)/t\Omega_{\rm Vir}(\lambda,0)$ is a one-dimensional trivial ${\rm Vir}$-module.
\item[{\rm(3)}] $\Omega_{\rm Vir}(\lambda_1,\alpha_1)\cong\Omega_{\rm Vir}(\lambda_2,\alpha_2)$ if and only if $\lambda_1=\lambda_2$ and $\alpha_1=\alpha_2$.
\end{itemize}
\end{lemm}

\subsection{An elementary result}

\begin{lemm}\label{Sequence-problem}
Let $\{\beta_m\,|\,m\in\Z_+\}$ be a sequence of complex numbers satisfying $\beta_0 = 1$ and
\begin{equation}\label{Sequence-problem-1}
\beta_{m+n} + (n-m)\beta_{m+n+1} = \beta_n\beta_m + n\beta_m\beta_{n+1} - m\beta_n\beta_{m+1}.
\end{equation}
Then $\beta_m = \beta^m$, where $\beta=\beta_1$.
\end{lemm}

\ni{\it Proof.}\ \
We prove this lemma by induction on $m$.
First, the cases for $m=0,1$ are clear.
Taking $(m,n)=(1,1)$ in \eqref{Sequence-problem-1}, we immediately see that
$\beta_{2} = \beta_1^2$.
Then, taking $(m,n)=(1,2)$ and $(2,1)$ respectively in \eqref{Sequence-problem-1}, we have
\begin{eqnarray*}
\label{Sequence-m3-1} \beta_{3} + \beta_{4} &\!\!\!=\!\!\!& \beta_1^3 + 2\beta_1\beta_3 - \beta_2^2, \\
\label{Sequence-m3-2} \beta_{3} - \beta_{4} &\!\!\!=\!\!\!& \beta_1^3 + \beta_2^2 - 2\beta_1\beta_3.
\end{eqnarray*}
Adding the above two equations, we see that $\beta_{3} = \beta_1^3$.

Let $m\ge 4$. Assume that the lemma holds for $k\le m-1$. Namely, we assume $\beta_k = \beta_1^k$ for $k\le m-1$.
Next, we consider the case for $m$.
Taking $(m,n)\rightarrow(m-2,1)$ in \eqref{Sequence-problem-1}, we have
\begin{equation*}
\beta_{m-1} + (3-m)\beta_{m} = \beta_1\beta_{m-2} + \beta_{m-2}\beta_2 - (m-2)\beta_1\beta_{m-1}.
\end{equation*}
A direct computation shows that $\beta_{m}= \beta_1^{m}$.
This completes the proof.
\QED

\section{\large{Free $\mathcal{U}(\mathfrak h)$-modules of rank one over $\mathcal{W}(1)$}}

Note that the Lie algebra $\mathcal{W}(1)$ is $\Z\times\Z_+$-graded, and recall that $\mathfrak h=\C L_{0,0}$ is the $(0,0)$-graded part.
In this section, we completely classify the free $\mathcal{U}(\mathfrak h)$-modules of rank one over $\mathcal{W}(1)$.

\subsection{Construction of free $\mathcal{U}(\mathfrak h)$-modules of rank one over $\mathcal{W}(1)$}

Let $\lambda\in\C^*$ and $\alpha, \beta\in\C$. Define the action of $\mathcal{W}(1)$ on the vector space of polynomials in one variable
$$
\Omega_{\mathcal{W}(1)}(\lambda, \alpha, \beta):=\C[t]
$$
by
\begin{equation}\label{action-of-w1}
L_{i,m}\cdot t^k = \sum_{s=0}^{{\rm min}\{m,\,k\}}s!\binom{m}{s}\binom{k}{s}\lambda^i\beta^{m-s-1}((m-s)\alpha - i\alpha\beta + \beta t)(t-i)^{k-s},
\end{equation}
where $i\in\Z, m\in\Z_+$ and $t^k\in\C[t]$. It is important to note here that we allow $\beta=0$:
if $\beta=0$ and $s=m$ in the sum of \eqref{action-of-w1}, although $\beta^{m-s-1}=\beta^{-1}$ formally appears, we view
$$
\beta^{m-s-1}((m-s)\alpha - i\alpha\beta + \beta t) =t - i\alpha
$$
as a whole, and thus \eqref{action-of-w1} still make sense.
Note further that if we take $i=m=0$ in \eqref{action-of-w1}, then we have $L_{0,0}\cdot t^k=t^{k+1}$.
Hence, $\Omega_{\mathcal{W}(1)}(\lambda, \alpha, \beta)$ is a $\mathcal{U}(\mathfrak h)$-module, which is free of rank one.
Furthermore, we find that $\Omega_{\mathcal{W}(1)}(\lambda, \alpha, \beta)$ is in fact a $\mathcal{W}(1)$-module.

\begin{prop}\label{W1-module}
The space $\Omega_{\mathcal{W}(1)}(\lambda, \alpha, \beta)$ is a $\mathcal{W}(1)$-module under the action \eqref{action-of-w1}.
\end{prop}

\ni{\it Proof.}\ \
By \eqref{brackets-for-w}, we see that $\{L_{i,0},L_{i,1}\mid i\in\Z\}$ is a generating set of $\mathcal{W}(1)$.
Hence, it is enough to check the commutation relations involving $L_{i,0}$ and $L_{i,1}$ on $t^k$ for any $k\in\Z_+$.
For $m=0,1,2$, the formula \eqref{action-of-w1} can be written explicitly as
\begin{eqnarray}
\label{action-i0-tk} L_{i,0}\cdot t^k &\!\!\!=\!\!\!& \lambda^i(t-i\alpha)(t-i)^k,\\
\label{action-i1-tk} L_{i,1}\cdot t^k &\!\!\!=\!\!\!&
\left\{
\begin{aligned}
&\lambda^i(\alpha-i\alpha\beta+\beta t), & k=0,\\
&\lambda^i(\alpha-i\alpha\beta+\beta t)(t-i)^k + k\lambda^i(t-i\alpha)(t-i)^{k-1}, & k\ge 1,
\end{aligned}
\right.\\
\label{action-i2-tk} L_{i,2}\cdot t^k &\!\!\!=\!\!\!&
\left\{
\begin{aligned}
&\lambda^i \beta(2\alpha-i\alpha\beta+\beta t), & k=0,\\
&\lambda^i \beta(2\alpha-i\alpha\beta+\beta t)(t-i) + 2\lambda^i(\alpha-i\alpha\beta+\beta t), & k=1,\\
& \begin{aligned}
&\lambda^i \beta(2\alpha-i\alpha\beta+\beta t)(t-i)^k + 2k\lambda^i(\alpha-i\alpha\beta+\beta t)(t-i)^{k-1} \\[-4pt]
& +k(k-1)\lambda^i(t-i\alpha)(t-i)^{k-2},\end{aligned} & k\ge 2.\\
\end{aligned}
\right.
\end{eqnarray}
Note that \eqref{action-i0-tk} is exactly the $\Vir$-module structure defined by \eqref{def-vir-mod-simple}.
Hence, by \cite{TZ2015}, we always have  $[L_{i,0}, L_{j,0}]\cdot t^k = L_{i,0}\cdot L_{j,0}\cdot t^k - L_{j,0}\cdot L_{i,0}\cdot t^k$.

Next, we prove $[L_{i,0},L_{j,1}]\cdot t^k = L_{i,0}\cdot L_{j,1}\cdot t^k - L_{j,1}\cdot L_{i,0}\cdot t^k$.
We only consider the case $k \geq 1$; the case $k=0$ can be checked more easily (the details are omitted).
On the one hand, by \eqref{action-i0-tk} and \eqref{action-i1-tk}, we have
\begin{eqnarray*}
\frac{1}{\lambda^{i+j}}[L_{i,0},L_{j,1}]\cdot t^k &\!\!\!=\!\!\!& \frac{1}{\lambda^{i+j}}(j-i)L_{i+j,1}\cdot t^k - \frac{1}{\lambda^{i+j}} L_{i+j,0}\cdot t^k \\
&\!\!\!=\!\!\!& (j-i)(\alpha-(i+j)\alpha\beta + \beta t)(t-i-j)^k \\
&& + k(j-i)(t-(i+j)\alpha)(t-i-j)^{k-1}\\
&& - (t-(i+j)\alpha)(t-i-j)^k.
\end{eqnarray*}
On the other hand, we have
\begin{eqnarray*}
\frac{1}{\lambda^{i+j}}L_{i,0}\cdot L_{j,1}\cdot t^k &\!\!\!=\!\!\!& \frac{1}{\lambda^{i}}L_{i,0}\cdot \left((\alpha-j\alpha\beta + \beta t)(t-j)^k + k(t-j\alpha)(t-j)^{k-1}\right)\\
&\!\!\!=\!\!\!& (t-i\alpha)(\alpha-j\alpha\beta + \beta (t-i))(t-i-j)^k\\
&& + k(t-i\alpha)(t-i-j\alpha)(t-i-j)^{k-1},
\end{eqnarray*}
and
\begin{eqnarray*}
\frac{1}{\lambda^{i+j}}L_{j,1}\cdot L_{i,0}\cdot t^k &\!\!\!=\!\!\!& \frac{1}{\lambda^{j}}L_{j,1}\cdot (t-i\alpha)(t-i)^k\\
&\!\!\!=\!\!\!& (\alpha-j\alpha\beta+\beta t)(t-j-i\alpha)(t-i-j)^k\\
&& + (t-j\alpha)(t-i-j)^k + k(t-j\alpha)(t-i\alpha-j)(t-i-j)^{k-1}.
\end{eqnarray*}
The above three formulas show that $[L_{i,0},L_{j,1}]\cdot t^k = L_{i,0}\cdot L_{j,1}\cdot t^k - L_{j,1}\cdot L_{i,0}\cdot t^k$.

Finally, we prove $[L_{i,1},L_{j,1}]\cdot t^k = L_{i,1}\cdot L_{j,1}\cdot t^k - L_{j,1}\cdot L_{i,1}\cdot t^k$.
As above, we only consider the most complicated case $k \geq 2$.
On the one hand, by \eqref{action-i2-tk}, we have
\begin{eqnarray*}
\frac{1}{\lambda^{i+j}}[L_{i,1},L_{j,1}]\cdot t^k &\!\!\!=\!\!\!& \frac{1}{\lambda^{i+j}}(j-i)L_{i+j,2}\cdot t^k\\
&\!\!\!=\!\!\!& \beta(j-i)(2\alpha-(i+j)\alpha\beta+\beta t)(t-i-j)^k \\
&&+ 2k(j-i)(\alpha-(i+j)\alpha\beta+\beta t)(t-i-j)^{k-1} \\
&&+ k(k-1)(j-i)(t-(i+j)\alpha)(t-i-j)^{k-2}.
\end{eqnarray*}
On the other hand, by \eqref{action-i1-tk}, we have
\begin{eqnarray*}
\frac{1}{\lambda^{i+j}}L_{i,1}\cdot L_{j,1}\cdot t^k &\!\!\!=\!\!\!& \frac{1}{\lambda^{i}}L_{i,1}\cdot \left((\alpha-j\alpha\beta+\beta t)(t-j)^k + k(t-j\alpha)(t-j)^{k-1}\right)\\
&\!\!\!=\!\!\!& (\alpha\!-\!i\alpha\beta\!+\!\beta t)\Big((\alpha\!-\!j\alpha\beta+\beta (t\!-\!i))(t\!-\!i\!-\!j)^k + k(t\!-\!i\!-\!j\alpha)(t\!-\!i\!-\!j)^{k-1}\Big)\\
&& + (t-i\alpha)\Big(\beta(t-i-j)^k + k(\alpha-j\alpha\beta+\beta (t-i))(t-i-j)^{k-1}\\
&& + k(t-i-j)^{k-1} + k(k-1)(t-i-j\alpha)(t-i-j)^{k-2}\Big).
\end{eqnarray*}
Using the above two formulas, it is straightforward to check the required relation.
\QED

\subsection{Classification of free $\mathcal{U}(\mathfrak h)$-modules of rank one over $\mathcal{W}(1)$}

Now we show that $\Omega_{\mathcal{W}(1)}(\lambda, \alpha, \beta)$ constructed in \eqref{action-of-w1} in fact exhaust all free $\mathcal{U}(\mathfrak h)$-modules of rank one over $\mathcal{W}(1)$. This is our first main result.

\begin{theo}\label{thm-1}
Any free $\mathcal{U}(\mathfrak h)$-module of rank one over $\mathcal{W}(1)$ is isomorphic to $\Omega_{\mathcal{W}(1)}(\lambda, \alpha, \beta)$
defined by \eqref{action-of-w1} for some $\lambda \in\C^*$ and $\alpha, \beta\in\C$.
\end{theo}

\ni{\it Proof.}\ \
Let $M$ be a free $\mathcal{U}(\mathfrak h)$-module of rank one over $\mathcal{W}(1)$.
By viewing $M$ as a $\Vir$-module, from Lemma~\ref{lemm-12}, we may assume that $M=\C[t]$ and
there exist some $\lambda\in\C^*$ and $\alpha\in\C$ such that the action of $\Vir\subseteq\mathcal{W}(1)$ on $M$ is as follows
\begin{equation}\label{action-i0}
L_{i,0}\cdot t^k=\lambda^i(t-i\alpha)(t-i)^k.
\end{equation}
We reduce the remaining proof to the following Lemma~\ref{w1-action-on-tk} and Lemma~\ref{w1-action-on-1}.
\QED

\begin{lemm}\label{w1-action-on-tk}
The action of $L_{i,m}$ on $t^k$ is a linear combination of the actions of $L_{i,j}$ ($m-{\rm min}\{m,k\} \leq j \leq m$) on $1$,
and more precisely we have
\begin{equation}\label{claim-relation}
L_{i,m}\cdot t^k = \sum_{s=0}^{{\rm min}\{m,\,k\}}s!\binom{m}{s}\binom{k}{s}(t-i)^{k-s} L_{i,m-s}\cdot 1.
\end{equation}
\end{lemm}
\ni{\it Proof.}\ \
We shall prove \eqref{claim-relation} by induction on $k$.
First, by \eqref{brackets-for-w}, we have
\begin{equation*}
L_{i,m}\cdot t = L_{i,m}\cdot L_{0,0}\cdot 1= (t-i)L_{i,m}\cdot 1 + mL_{i,m-1}\cdot 1,
\end{equation*}
which implies that \eqref{claim-relation} holds for $k=1$.
Let $k\ge 2$. Assume that \eqref{claim-relation} holds for $k-1$, namely we have
\begin{equation}\label{case-k-1}
L_{i,m}\cdot t^{k-1} = \sum_{s=0}^{{\rm min}\{m,\,k-1\}}s!\binom{m}{s}\binom{k-1}{s}(t-i)^{k-s-1} L_{i,m-s}\cdot 1.
\end{equation}
Next, we prove that \eqref{claim-relation} holds for $k$ in two cases.

{\bf Case 1:} $m<k$. In this case, by \eqref{case-k-1}, we have
\begin{eqnarray*}
L_{i,m}\cdot t^k &\!\!\!=&\!\!\! L_{i,m}\cdot L_{0,0}\cdot t^{k-1}=(t-i)L_{i,m}\cdot t^{k-1} + mL_{i,m-1}\cdot t^{k-1}\\
&\!\!\!=&\!\!\! (t-i)\sum_{s=0}^{{\rm min}\{m,\,k-1\}}s!\binom{m}{s}\binom{k-1}{s}(t-i)^{k-s-1} L_{i,m-s}\cdot 1\\
&\!\!\!&\!\!\! + m\sum_{s=0}^{{\rm min}\{m-1,\,k-1\}}s!\binom{m-1}{s}\binom{k-1}{s}(t-i)^{k-s-1} L_{i,m-s-1}\cdot 1\\
&\!\!\!=&\!\!\! (t-i)^{k}L_{i,m}\cdot 1 + \sum_{s=1}^{m}s!\binom{m}{s}\binom{k-1}{s}(t-i)^{k-s} L_{i,m-s}\cdot 1\\
&\!\!\!&\!\!\! + m\sum_{s=0}^{m-1}s!\binom{m-1}{s}\binom{k-1}{s}(t-i)^{k-s-1} L_{i,m-s-1}\cdot 1\\
&\!\!\!=&\!\!\! (t-i)^{k}L_{i,m}\cdot 1 + \sum_{s=1}^{m}s!\binom{m}{s}\binom{k-1}{s}(t-i)^{k-s} L_{i,m-s}\cdot 1\\
&\!\!\!&\!\!\! + \sum_{s=1}^{m}s!\binom{m}{s}\binom{k-1}{s-1}(t-i)^{k-s} L_{i,m-s}\cdot 1\\
&\!\!\!=&\!\!\! \sum_{s=0}^{m}s!\binom{m}{s}\binom{k}{s}(t-i)^{k-s} L_{i,m-s}\cdot 1.
\end{eqnarray*}
Here we have used the combinatoric formula (cf.~\eqref{formula})
\begin{equation*}
\binom{k-1}{s} + \binom{k-1}{s-1} = \binom{k}{s}
\end{equation*}
in the last equality.

{\bf Case 2:} $m\geq k$.
Similarly, in this case, we have
\begin{eqnarray*}
L_{i,m}\cdot t^k &\!\!\!=&\!\!\! L_{i,m}\cdot L_{0,0}\cdot t^{k-1}=(t-i)L_{i,m}\cdot t^{k-1} + mL_{i,m-1}\cdot t^{k-1}\\
&\!\!\!=&\!\!\! (t-i)\sum_{s=0}^{{\rm min}\{m,\,k-1\}}s!\binom{m}{s}\binom{k-1}{s}(t-i)^{k-s-1} L_{i,m-s}\cdot 1\\
&\!\!\!&\!\!\! + m\sum_{s=0}^{{\rm min}\{m-1,\,k-1\}}s!\binom{m-1}{s}\binom{k-1}{s}(t-i)^{k-s-1} L_{i,m-s-1}\cdot 1\\
&\!\!\!=&\!\!\! (t-i)^{k}L_{i,m}\cdot 1 + \sum_{s=1}^{k-1}s!\binom{m}{s}\binom{k-1}{s}(t-i)^{k-s} L_{i,m-s}\cdot 1\\
&\!\!\!&\!\!\! + m\sum_{s=0}^{k-2}s!\binom{m-1}{s}\binom{k-1}{s}(t-i)^{k-s-1} L_{i,m-s-1}\cdot 1 + m(k-1)!\binom{m-1}{k-1}L_{i,m-k}\cdot 1\\
&\!\!\!=&\!\!\! (t-i)^{k}L_{i,m}\cdot 1 + \sum_{s=1}^{k-1}s!\binom{m}{s}\binom{k-1}{s}(t-i)^{k-s} L_{i,m-s}\cdot 1\\
&\!\!\!&\!\!\! + \sum_{s=1}^{k-1}s!\binom{m}{s}\binom{k-1}{s-1}(t-i)^{k-s} L_{i,m-s}\cdot 1 + m(k-1)!\binom{m-1}{k-1}L_{i,m-k}\cdot 1\\
&\!\!\!=&\!\!\! \sum_{s=0}^{k}s!\binom{m}{s}\binom{k}{s}(t-i)^{k-s} L_{i,m-s}\cdot 1.
\end{eqnarray*}
This completes the proof.
\QED

\begin{lemm}\label{w1-action-on-1}
There exists some $\beta\in\C$ such that $L_{i,m}\cdot 1 = \lambda^i\beta^{m-1}(m\alpha - i\alpha\beta + \beta t)$.
\end{lemm}
\ni{\it Proof.}\ \
Denote $F_{i,m}(t)= L_{i,m}\cdot 1$. We determine $F_{i,m}(t)$ by induction on $m$.

The case for $m=0$ has been given by \eqref{action-i0} with $k=0$
(note that, as before, if $m=\beta=0$ in this lemma, we have viewed $\beta^{-1}(m\alpha- i\alpha\beta+\beta t)=t- i\alpha$ as a whole).

Let us first determine $F_{i,1}(t)$. Applying
$$
[L_{0,1},L_{i,0}]=iL_{i,1}+L_{i,0}, \quad [L_{-i,0},L_{i,1}]=2iL_{0,1}-L_{0,0}
$$
respectively on $1$, by \eqref{action-i0}, we obtain
\begin{eqnarray}
\label{w+4.3}
iF_{i,1}(t) &\!\!\!=&\!\!\! \lambda^i(t-i\alpha)F_{0,1}(t)-\lambda^i(t-i\alpha)F_{0,1}(t-i) + \lambda^ii\alpha, \\
\label{w+4.4}
2iF_{0,1}(t) &\!\!\!=&\!\!\! \lambda^{-i}(t+i\alpha)F_{i,1}(t+i)-\lambda^{-i}(t+i\alpha-i)F_{i,1}(t)+i\alpha.
\end{eqnarray}
Multiplying \eqref{w+4.4} by $i$, and then using the relation \eqref{w+4.3}, we can derive that
\begin{eqnarray}
\nonumber  2i^2\alpha &\!\!\!=\!\!\!& 2(t^2-i^2\alpha^2+i^2\alpha+i^2)F_{0,1}(t)-(t^2+it-i^2\alpha^2+i^2\alpha)F_{0,1}(t+i)\\
\label{w+4.22}  &\!\!\!\!\!\!& -(t^2-it-i^2\alpha^2+i^2\alpha)F_{0,1}(t-i).
\end{eqnarray}
If $F_{0,1}(t)=0$, then by \eqref{w+4.3} and \eqref{w+4.22}, we see that $F_{i,1}(t)=0$ for all $i\in\Z$.
This proves this lemma for the case $(\alpha,\beta)=(0,0)$.
If $F_{0,1}(t)\ne 0$, we let $\deg\,F_{0,1}(t)=K$, and assume that
\begin{equation}\label{w+4.23}
F_{0,1}(t)=\sum\limits_{r=0}^{K}a_{r}t^{r},\quad \mbox{where} \ \ a_r\in\C \ \mbox{and}\ a_K\neq0.
\end{equation}
If $K=0$, then $F_{0,1}(t)=a_0$. By \eqref{w+4.22} with $i=1$, we have $a_0=\alpha$, and thus $F_{0,1}(t)=\alpha$ ($\ne 0$).
If $K=1$, then $F_{0,1}(t)=a_0+a_1 t$. By \eqref{w+4.22}, one can also derive that $a_0=\alpha$. Redenote $a_1$ by $\beta$, we have $F_{0,1}(t)=\alpha+\beta t$.
If $K\ge 2$, substituting \eqref{w+4.23} into \eqref{w+4.22} with $i=1$ and comparing the coefficients of $t^K$ on both sides, we obtain
$(K^2 + K - 2) a_K=0$, a contradiction.
Hence, in general, we have $F_{0,1}(t) = \alpha + \beta t$. Then, by \eqref{w+4.3}, we obtain
$F_{i,1}(t) = \lambda^i(\alpha - i\alpha\beta + \beta t)$.
This proves this lemma for the case $(\alpha,\beta)\ne (0,0)$.

Let $m\ge 1$. Assume that this lemma holds for $0\leq k\leq m-1$, namely,
\begin{equation}\label{case-k}
F_{i,k}(t) = \lambda^i\beta^{k-1}(k\alpha - i\alpha\beta + \beta t), \  0\leq k\leq m-1.
\end{equation}
Next, we determine $F_{i,m}(t)$.
Using the same arguments as for the case $F_{i,1}(t)$, by relations
$$
[L_{0,m},L_{i,0}] = iL_{i,m} + m L_{i,m-1},\quad  [L_{-i,0},L_{i,m}] = 2iL_{0,m} - m L_{0,m-1}
$$
and \eqref{case-k}, one can derive that
$F_{i,m}(t)=\lambda^i (m\alpha\beta^{m-1} - i\alpha b_{m} + b_{m}t)$ for some $b_{m}\in\C$.
Then, applying
$$
[L_{1,1}, L_{0,m-1}]=-L_{1,m}- (m-2)L_{1,m-1}
$$
on $1$, we obtain $b_{m} = \beta^m$, and thus
$F_{i,m}(t) = \lambda^i\beta^{m - 1}(m\alpha - i\alpha\beta + \beta t)$. This completes the proof.
\QED

\subsection{Simplicity and isomorphism classification of $\mathcal{W}(1)$-modules}

Next, we determine the simplicity of $\Omega_{\mathcal{W}(1)}(\lambda,\alpha,\beta)$.
One will see that although the results are similar to Lemmas~\ref{lemm-13}(1) and (2), the proofs, especially for Theorem~\ref{w1-simplicity}(2), are rather non-trivial.

\begin{theo}\label{w1-simplicity}
\parskip-3pt
\begin{itemize}\parskip-3pt
\item[{\rm(1)}] $\Omega_{\mathcal{W}(1)}(\lambda,\alpha,\beta)$ is simple if and only if $\alpha\neq0$.
\item[{\rm(2)}] $\Omega_{\mathcal{W}(1)}(\lambda,0,\beta)$ has a unique proper submodule $t\Omega_{\mathcal{W}(1)}(\lambda,0,\beta)\cong\Omega_{\mathcal{W}(1)}(\lambda,1,\beta)$, and the quotient $\Omega_{\mathcal{W}(1)}(\lambda,0,\beta)/t\Omega_{\mathcal{W}(1)}(\lambda,0,\beta)$ is a one-dimensional trivial ${\mathcal{W}(1)}$-module.
\end{itemize}
\end{theo}

\ni{\it Proof.}\ \
(1) Let $M=\Omega_{\mathcal{W}(1)}(\lambda,\alpha,\beta)$. If $\alpha\neq0$, by viewing $M$ as a ${\rm Vir}$-module, from Lemma~\ref{lemm-13}(1), we see that $M$ is simple.
If $\alpha=0$, one can easily see that $t M$ is a submodule of $M$.

(2) From definition, one can easily see that $t\Omega_{\mathcal{W}(1)}(\lambda,0,\beta)$ is the unique proper submodule of $\Omega_{\mathcal{W}(1)}(\lambda,0,\beta)$.
Next, we prove $t\Omega_{\mathcal{W}(1)}(\lambda,0,\beta)\cong\Omega_{\mathcal{W}(1)}(\lambda,1,\beta)$ by comparing the actions of $\mathcal{W}(1)$ on these two modules.

First, we consider the action of $\mathcal{W}(1)$ on $t\Omega_{\mathcal{W}(1)}(\lambda,0,\beta)$.
For $k\ge 0$, by \eqref{action-of-w1}, we have
\begin{equation}\label{w1-on-t0}
L_{i,m}\cdot (t\cdot t^k) = t\left(\sum_{s=0}^{{\rm min}\{m,\,k+1\}}s!\binom{m}{s}\binom{k+1}{s}\lambda^i\beta^{m-s}(t-i)^{k-s+1}\right).
\end{equation}
Next, we consider the action of $\mathcal{W}(1)$ on $\Omega_{\mathcal{W}(1)}(\lambda,1,\beta)$ in two cases. Let $k\ge 0$.

{\bf Case 1:} $m\le k$. In this case, by \eqref{action-of-w1}, we have
\begin{eqnarray*}
L_{i,m}\cdot t^k &\!\!\!=\!\!\!& \sum_{s=0}^{m}s!\binom{m}{s}\binom{k}{s}\lambda^i\beta^{m-s-1}(\beta (t-i) + (m-s))(t-i)^{k-s}\\
&\!\!\!=\!\!\!& \sum_{s=0}^{m}s!\binom{m}{s}\binom{k}{s}\lambda^i\beta^{m-s}(t-i)^{k-s+1} + \sum_{s=0}^{m}s!\binom{m}{s}\binom{k}{s}\lambda^i\beta^{m-s-1}(m-s)(t-i)^{k-s}\\
&\!\!\!=\!\!\!&  \lambda^i\beta^{m}(t-i)^{k+1} + \sum_{s=1}^{m}s!\binom{m}{s}\binom{k}{s}\lambda^i\beta^{m-s}(t-i)^{k-s+1}\\
&\!\!\!\!\!\!& +  \sum_{s=0}^{m-1}s!\binom{m}{s}\binom{k}{s}\lambda^i\beta^{m-s-1}(m-s)(t-i)^{k-s}\\
&\!\!\!=\!\!\!&  \lambda^i\beta^{m}(t-i)^{k+1} + \sum_{s=1}^{m}s!\binom{m}{s}\binom{k}{s}\lambda^i\beta^{m-s}(t-i)^{k-s+1}\\
&\!\!\!\!\!\!& +  \sum_{s=1}^{m}(s-1)!\binom{m}{s-1}\binom{k}{s-1}\lambda^i\beta^{m-s}(m-s+1)(t-i)^{k-s+1}\\
&\!\!\!=\!\!\!&  \lambda^i\beta^{m}(t-i)^{k+1} + \sum_{s=1}^{m}s!\binom{m}{s}\binom{k}{s}\lambda^i\beta^{m-s}(t-i)^{k-s+1}\\
&\!\!\!\!\!\!& + \sum_{s=1}^{m}s!\binom{m}{s}\binom{k}{s-1}\lambda^i\beta^{m-s}(t-i)^{k-s+1}\\
&\!\!\!=\!\!\!&\sum_{s=0}^{m}s!\binom{m}{s}\binom{k+1}{s}\lambda^i\beta^{m-s}(t-i)^{k-s+1}.
\end{eqnarray*}
Here, we have again used the combinatoric formula \eqref{formula} in the last equality.
\vskip5pt

{\bf Case 2:} $m> k$.
Similarly, in this case, we have
\begin{eqnarray*}
L_{i,m}\cdot t^k &\!\!\!=\!\!\!& \sum_{s=0}^{k}s!\binom{m}{s}\binom{k}{s}\lambda^i\beta^{m-s-1}(\beta (t-i) + (m-s))(t-i)^{k-s}\\
&\!\!\!=\!\!\!& \sum_{s=0}^{k}s!\binom{m}{s}\binom{k}{s}\lambda^i\beta^{m-s}(t-i)^{k-s+1} +\sum_{s=0}^{k}s!\binom{m}{s}\binom{k}{s}\lambda^i\beta^{m-s-1}(m-s)(t-i)^{k-s}\\
&\!\!\!=\!\!\!&  \lambda^i\beta^{m}(t-i)^{k+1} + \sum_{s=1}^{k}s!\binom{m}{s}\binom{k}{s}\lambda^i\beta^{m-s}(t-i)^{k-s+1}\\
&\!\!\!\!\!\!& +  \sum_{s=0}^{k-1}s!\binom{m}{s}\binom{k}{s}\lambda^i\beta^{m-s-1}(m-s)(t-i)^{k-s} +k!\binom{m}{k}\lambda^i\beta^{m-k-1}(m-k)\\
&\!\!\!=\!\!\!&  \lambda^i\beta^{m}(t-i)^{k+1} + \sum_{s=1}^{k}s!\binom{m}{s}\binom{k}{s}\lambda^i\beta^{m-s}(t-i)^{k-s+1}\\
&\!\!\!\!\!\!& +  \sum_{s=1}^{k}(s-1)!\binom{m}{s-1}\binom{k}{s-1}\lambda^i\beta^{m-s}(m-s+1)(t-i)^{k-s+1}\\
&\!\!\!\!\!\!&  + (k+1)!\binom{m}{k+1}\lambda^i\beta^{m-k-1}\\
&\!\!\!=\!\!\!&  \lambda^i\beta^{m}(t-i)^{k+1} + \sum_{s=1}^{k}s!\binom{m}{s}\binom{k}{s}\lambda^i\beta^{m-s}(t-i)^{k-s+1}\\
&\!\!\!\!\!\!& + \sum_{s=1}^{k}s!\binom{m}{s}\binom{k}{s-1}\lambda^i\beta^{m-s}(t-i)^{k-s+1} + (k+1)!\binom{m}{k+1}\lambda^i\beta^{m-k-1}\\
&\!\!\!=\!\!\!&\sum_{s=0}^{k+1}s!\binom{m}{s}\binom{k+1}{s}\lambda^i\beta^{m-s}(t-i)^{k-s+1}.
\end{eqnarray*}
Summarizing the above two cases, we have
\begin{equation}\label{w1-on-1}
L_{i,m}\cdot t^k = \sum_{s=0}^{{\rm min}\{m,\,k+1\}}s!\binom{m}{s}\binom{k+1}{s}\lambda^i\beta^{m-s}(t-i)^{k-s+1}.
\end{equation}
Comparing \eqref{w1-on-t0} with \eqref{w1-on-1}, we see that $t\Omega_{\mathcal{W}(1)}(\lambda,0,\beta)\cong\Omega_{\mathcal{W}(1)}(\lambda,1,\beta)$.

At last, it is clear that the quotient $\Omega_{\mathcal{W}(1)}(\lambda,0,\beta)/t\Omega_{\mathcal{W}(1)}(\lambda,0,\beta)$ is a one-dimensional trivial ${\mathcal{W}(1)}$-module.
This completes the proof.
\QED
\vskip15pt

The isomorphism classification of $\Omega_{\mathcal{W}(1)}(\lambda,\alpha,\beta)$ is as follows.

\begin{theo}\label{w1-iso-classification}
$\Omega_{\mathcal{W}(1)}(\lambda_1,\alpha_1,\beta_1)\cong\Omega_{\mathcal{W}(1)}(\lambda_2,\alpha_2,\beta_2)$
if and only if $\lambda_1=\lambda_2$, $\alpha_1=\alpha_2$ and $\beta_1=\beta_2$.
\end{theo}

\ni{\it Proof.}\ \
Suppose that $\varphi$ is a module isomorphism from $\Omega_{\mathcal{W}(1)}(\lambda_1,\alpha_1,\beta_1)$ to $\Omega_{\mathcal{W}(1)}(\lambda_2,\alpha_2,\beta_2)$.
Denote by $\varphi^{-1}$ the inverse of $\varphi$. Let $f(t)=\varphi^{-1}(1)$. Since $L_{0,0}\cdot t^k=t^{k+1}$, we have
$$
1=\varphi(f(t))=\varphi(f(L_{0,0})\cdot1)=f(L_{0,0})\cdot\varphi(1)=f(t)\varphi(1).
$$
Hence, $\varphi(1)\in\C^*$. Computing $\varphi(L_{i,0}\cdot1)$, we have
$$
\varphi(L_{i,0}\cdot1)=\varphi(\lambda_1^i(t-i\alpha_1))=\varphi(\lambda_1^i(L_{0,0}-i\alpha_1)\cdot1)=\lambda_1^i(L_{0,0}-i\alpha_1)\cdot\varphi(1)=\lambda_1^i(t-i\alpha_1)\varphi(1).
$$
On the other hand, we have (note that $\varphi(1)\in\C^*$)
$$
\varphi(L_{i,0}\cdot1)=L_{i,0}\cdot\varphi(1)=\lambda_2^i(t-i\alpha_2)\varphi(1).
$$
Comparing the above two formulas, we must have $\lambda_1=\lambda_2$ and $\alpha_1=\alpha_2$.
Similarly, computing $\varphi(L_{i,1}\cdot1)$, one can derive that
\begin{equation*}
\lambda_1^i(\alpha_1 - i\alpha_1\beta_1 + \beta_1 t) \varphi(1) = \lambda_2^i(\alpha_2 - i\alpha_2\beta_2 + \beta_2 t) \varphi(1),
\end{equation*}
which implies that $\beta_1=\beta_2$.
This completes the proof.
\QED

\section{\large{Free $\mathcal{U}(\mathfrak h)$-modules of rank one over $\mathcal{W}(-1)$}}

Recall that $\mathfrak h=\C L_{0,0}$. In this section, we completely classify the free $\mathcal{U}(\mathfrak h)$-modules of rank one over $\mathcal{W}(-1)$.
Comparing with Section~3, one will see that the problem becomes more difficult,
especially when determining the actions on $1$ (one can compare the following Lemma~\ref{w-1-action-on-1} with Lemma~\ref{w1-action-on-1}).
Thus we need more technical analysis.

\subsection{\large{Construction of free $\mathcal{U}(\mathfrak h)$-modules of rank one over $\mathcal{W}(-1)$}}

Let $\lambda\in\C^*$ and $\alpha,\beta\in\C$.
We define the action of $\mathcal{W}(-1)$ on the vector space of polynomials in one variable
$$
\Omega_{\mathcal{W}(-1)}(\lambda,\alpha,\beta):=\C[t]
$$
by
\begin{eqnarray}
\begin{aligned}\label{action-of-w-1}
L_{i,m}\cdot t^k &= \sum_{s=0}^{k}(-1)^ss!\binom{m+s-1}{s}\binom{k}{s}\lambda^i\beta^{m+s}(t-i\alpha-(m+s)\alpha\beta)(t-i)^{k-s},
\end{aligned}
\end{eqnarray}
where $i\in\Z$, $m\in\Z_+$ and $t^k\in\C[t]$. It should note here that $\binom{-1}{0} = 1$ and $\binom{n-1}{n} = 0$ if $n> 0$.
This guarantees that $\Omega_{\mathcal{W}(-1)}(\lambda, \alpha, \beta)$ is a $\mathcal{U}(\mathfrak h)$-module, which is free of rank one.
Similar to Proposition~\ref{W1-module}, one can further prove that $\Omega_{\mathcal{W}(-1)}(\lambda,\alpha,\beta)$
is a $\mathcal{W}(-1)$-module.

\begin{prop}\label{W-1-module}
The space $\Omega_{\mathcal{W}(-1)}(\lambda, \alpha, \beta)$ is a $\mathcal{W}(-1)$-module under the action \eqref{action-of-w-1}.
\end{prop}

\ni{\it Proof.}\ \
Note first that \eqref{action-of-w-1} is equivalent to
\begin{equation}\label{action-of-w-1-ft}
L_{i,m}\cdot f(t) = \sum_{s=0}^{+\infty}(-1)^s\binom{m+s-1}{s}\lambda^i\beta^{m+s}(t-i\alpha-(m+s)\alpha\beta)f^{(s)}(t-i),
\end{equation}
where $f(t)\in\C[t]$ and $f^{(n)}(t)$ stands for the $n$-order derivative of $f(t)$.
Similar to case for $\mathcal{W}(1)$, the set $\{L_{i,0},L_{i,1}\mid i\in\Z\}$ is also a generating set of $\mathcal{W}(-1)$,
and thus it is enough to check the commutation relations involving $L_{i,0}$ and $L_{i,1}$ on $f(t)$.
Here we only check the most complicated case: the commutation relations between $L_{i,1}$'s on $f(t)$
(we leave the details of other cases to interested readers).
On the one hand, by \eqref{action-of-w-1-ft} with $m=2$, we have
\begin{eqnarray*}
\frac{1}{\lambda^{i+j}}[L_{i,1},L_{j,1}]\cdot f(t) = \frac{1}{\lambda^{i+j}}(j-i)L_{i+j,2}\cdot f(t) = t X_{L}(t) + \alpha Y_{L}(t),
\end{eqnarray*}
where
\begin{eqnarray*}
X_{L}(t) &\!\!\!=\!\!\!& (j-i)\sum_{s=0}^{+\infty}(-1)^s(s+1)\beta^{s+2}f^{(s)}(t-i-j),\\
Y_{L}(t) &\!\!\!=\!\!\!& (j-i)\sum_{s=0}^{+\infty}(-1)^{s+1}(s+1)\beta^{s+2}(i+j+(s+2)\beta)f^{(s)}(t-i-j).
\end{eqnarray*}
On the other hand, by \eqref{action-of-w-1-ft} with $m=1$, we have
\begin{eqnarray*}
\frac{1}{\lambda^{i+j}}L_{i,1}\cdot L_{j,1}\cdot f(t)
&\!\!\!=\!\!\!& \frac{1}{\lambda^{i}}L_{i,1}\cdot \sum_{s=0}^{+\infty}(-1)^s\beta^{s+1}(t-j\alpha-(s+1)\alpha\beta)f^{(s)}(t-j)\\
&\!\!\!=\!\!\!& \sum_{s=0}^{+\infty}\sum_{r=0}^{+\infty}(-1)^{s+r}\beta^{s+r+2}(t-i\alpha-(r+1)\alpha\beta)\\
&& \times\left((t-i-j\alpha-(s+1)\alpha\beta)f^{(s+r)}(t-i-j) + rf^{(s+r-1)}(t-i-j)\right),
\end{eqnarray*}
and thus
\begin{eqnarray*}
&& \frac{1}{\lambda^{i+j}}L_{i,1}\cdot L_{j,1}\cdot f(t) - \frac{1}{\lambda^{i+j}}L_{j,1}\cdot L_{i,1}\cdot f(t)\\
&\!\!\!=\!\!\!& (j-i)\sum_{s=0}^{+\infty}\sum_{r=0}^{+\infty}(-1)^{s+r}\beta^{s+r+2}\left(t-(i+j)\alpha-(r+1)\alpha\beta+(r-s)\alpha^2\beta\right)f^{(s+r)}(t-i-j)\\
&& +(j-i)\sum_{s=0}^{+\infty}\sum_{r=0}^{+\infty}(-1)^{s+r}\beta^{s+r+2}\alpha rf^{(s+r-1)}(t-i-j)\\
&\!\!\!=\!\!\!& t X_{R}(t) + \alpha Y_{R}(t) + \alpha^2 Z(t),
\end{eqnarray*}
where
\begin{eqnarray*}
X_{R}(t) &\!\!\!=\!\!\!& (j-i)\sum_{s=0}^{+\infty}\sum_{r=0}^{+\infty}(-1)^{s+r}\beta^{s+r+2}f^{(s+r)}(t-i-j),\\
Y_{R}(t) &\!\!\!=\!\!\!& (j-i)\sum_{s=0}^{+\infty}\sum_{r=0}^{+\infty}(-1)^{s+r+1}\beta^{s+r+2}(i+j+(r+1)\beta)f^{(s+r)}(t-i-j) \\
         && +(j-i)\sum_{s=0}^{+\infty}\sum_{r=0}^{+\infty}(-1)^{s+r}\beta^{s+r+2} rf^{(s+r-1)}(t-i-j),\\
Z(t) &\!\!\!=\!\!\!& (j-i)\sum_{s=0}^{+\infty}\sum_{r=0}^{+\infty}(-1)^{s+r}\beta^{s+r+3}(r-s) f^{(s+r)}(t-i-j).
\end{eqnarray*}
For any $k\in\Z_+$, we observe that
\begin{eqnarray}
\label{observe1}(-1)^{k}(k+1) &\!\!\!=\!\!\!& \sum_{s+r=k}(-1)^{s+r}, \\
\label{observe2}(-1)^{k+1}(k+1) &\!\!\!=\!\!\!& \sum_{s+r=k}(-1)^{s+r+1}, \\
\label{observe3}(-1)^{k+1}(k+1)(k+2) &\!\!\!=\!\!\!& \sum_{s+r=k}(-1)^{s+r+1}(r+1)+\sum_{s+r=k+1}(-1)^{s+r}r.
\end{eqnarray}
By \eqref{observe1}, we see that $X_{L}(t)=X_{R}(t)$.
By \eqref{observe2} and \eqref{observe3}, we see that $Y_{L}(t)=Y_{R}(t)$.
By the symmetry of $s$ and $r$, we have $Z(t)=0$.
Hence, $[L_{i,1},L_{j,1}]\cdot f(t) = L_{i,1}\cdot L_{j,1}\cdot f(t) - L_{j,1}\cdot L_{i,1}\cdot f(t)$.
\QED

\subsection{\large{Classification of free $\mathcal{U}(\mathfrak h)$-modules of rank one over $\mathcal{W}(-1)$}}

The following is our second main result, which states that $\Omega_{\mathcal{W}(-1)}(\lambda,\alpha,\beta)$ exhausts all free $\mathcal{U}(\mathfrak h)$-modules of rank one over $\mathcal{W}(-1)$.

\begin{theo}\label{thm-2}
Any free $\mathcal{U}(\mathfrak h)$-module of rank one over $\mathcal{W}(-1)$ is isomorphic to $\Omega_{\mathcal{W}(-1)}(\lambda,\alpha,\beta)$
defined by \eqref{action-of-w-1} for some $\lambda\in\C^*$ and $\alpha,\beta\in\C$.
\end{theo}

\ni{\it Proof.}\ \
Let $M$ be a free $\mathcal{U}(\mathfrak h)$-module of rank one over $\mathcal{W}(-1)$.
By viewing $M$ as a $\Vir$-module, from Lemma~\ref{lemm-12}, we may assume that $M=\C[t]$ and
there exist some $\lambda\in\C^*$ and $\alpha\in\C$ such that the action of $\Vir\subseteq\mathcal{W}(-1)$ on $M$ is as follows
\begin{equation}\label{w-1-action-i0}
L_{i,0}\cdot t^k=\lambda^i(t-i\alpha)(t-i)^k.
\end{equation}
We reduce the remaining proof to the following Lemma~\ref{w-1-action-on-tk} and Lemma~\ref{w-1-action-on-1}.
\QED

\begin{lemm}\label{w-1-action-on-tk}
The action of $L_{i,m}$ on $t^k$ is a linear combination of the actions of $L_{i,j}$ ($m \leq j \leq m+k$) on $1$,
and more precisely we have
\begin{equation}\label{claim-relation-2.1}
L_{i,m}\cdot t^k = \sum_{s=0}^{k}(-1)^ss!\binom{m+s-1}{s}\binom{k}{s}(t-i)^{k-s}L_{i,m+s}\cdot 1.
\end{equation}
\end{lemm}
\ni{\it Proof.}\ \
The proof is similar to that of Lemma~\ref{w1-action-on-tk} and is omitted.
\QED

\begin{lemm}\label{w-1-action-on-1}
There exists some $\beta\in\C$ such that $L_{i,m}\cdot 1 = \lambda^i\beta^{m}(t-i\alpha - m\alpha\beta)$.
\end{lemm}
\ni{\it Proof.}\ \
Denote $G_{i,m}(t)= L_{i,m}\cdot 1$.
Since the Lie structure of $\mathcal{W}(-1)$ is essentially different from that of $\mathcal{W}(1)$, one cannot prove this lemma by induction on $m$ as in Lemma~\ref{w1-action-on-1}.
In fact, the situation becomes more difficult, and thus we need more technical analysis.
To simplify the proof, we shall use the shifted notation $Y_{i,m}(t)= \lambda^{-i}G_{i,m}(t)$.

Let $N_i=\deg\,Y_{i,1}(t)$, $i\in\Z$.
Assume that
\begin{equation}\label{Yi1t}
Y_{i,1}(t)=\sum_{r=0}^{N_i}a_r^{(i)} t^r,\quad \mbox{where} \ \ a_r^{(i)}\in\C \ \mbox{and}\ a_{N_i}^{(i)}\neq0.
\end{equation}
Let $f(x_1,x_2,\ldots,x_n)\in \C[x_1,x_2,\ldots,x_n]$ be a polynomial in $n$ variables.
In the following, we shall use the notation $\mathrm{Coeff}_{f(x_1,x_2,\ldots,x_n)}x_r^{i}$ to denote the coefficient of $x_r^{i}$ in $f(x_1,x_2,\ldots,x_n)\in \C[x_r]$.
In particular, if $i=0$, we use it to denote the constant term.

\begin{clai}\label{cliam-4.1}
We have $N_i=N_0$ and $a_{N_i}^{(i)}=a_{N_0}^{(0)}$ for all $i\in\Z$.
We shall simply denote $N_0$ by $N$, and $a_{N_0}^{(0)}$ by $a_N$.
\end{clai}

Applying
\begin{equation}\label{start-equations}
[L_{j,1},L_{i-j,0}] = (i-2j)L_{i,1} - L_{i,2}, \quad [L_{-j,0},L_{j,1}]= 2jL_{0,1} + L_{0,2}
\end{equation}
respectively on $1$, by \eqref{w-1-action-i0}, we obtain
\begin{eqnarray}
\label{Lj1,Li-j0} (t\!-\!j-(i\!-\!j)\alpha)Y_{j,1}(t)\! -\! (t\!-(i\!-\!j)\alpha)Y_{j,1}(t-i+j)\! -\! Y_{j,2}(t) &\!\!\!=\!\!\!& (i-2j)Y_{i,1}(t) - Y_{i,2}(t),\\
\label{L-i0,Li0}  (t+j\alpha)Y_{j,1}(t+j) - (t-j+j\alpha)Y_{j,1}(t) + Y_{j,2}(t) &\!\!\!=\!\!\!& 2jY_{0,1}(t) + Y_{0,2}(t).
\end{eqnarray}
Taking $j = 0$ in \eqref{Lj1,Li-j0}, we have
\begin{eqnarray}
\begin{aligned}\label{L01,Li0}
(t-i\alpha)Y_{0,1}(t) - (t-i\alpha)Y_{0,1}(t-i) - Y_{0,2}(t) = iY_{i,1}(t) - Y_{i,2}(t).
\end{aligned}
\end{eqnarray}
Computing $\eqref{Lj1,Li-j0}+\eqref{L-i0,Li0}-\eqref{L01,Li0}$, we obtain
\begin{eqnarray}
\nonumber  &\!\!\!\!\!\!& (t+j\alpha)Y_{j,1}(t+j) - i\alpha Y_{j,1}(t) - (t-(i-j)\alpha)Y_{j,1}(t-i+j) + 2jY_{i,1}(t)\\
\label{add-eqa}   &\!\!\!=\!\!\!& (t-i\alpha+2j)Y_{0,1}(t) - (t-i\alpha)Y_{0,1}(t-i).
\end{eqnarray}
Taking $j=i$ in \eqref{add-eqa}, we have
\begin{equation}\label{add-eqa-j=i}
(t+i\alpha)Y_{i,1}(t+i)  - (t+i\alpha-2i) Y_{i,1}(t) =(t-i\alpha+2i)Y_{0,1}(t) - (t-i\alpha)Y_{0,1}(t-i).
\end{equation}
Let ${\bf{L}}_1(t)$ and ${\bf{R}}_1(t)$ denote respectively the left- and right-hand sides of \eqref{add-eqa-j=i}.
Considering the coefficients of the highest degree terms, we have
$$
{\rm{Coeff}}_{{\bf{L}}_1(t)} t^{N_i}=(N_i+2) i a_{N_i}^{(i)}, \quad {\rm{Coeff}}_{{\bf{R}}_1(t)} t^{N_0}=(N_0+2) i a_{N_0}^{(0)}.
$$
By \eqref{add-eqa-j=i} with $i\ne 0$, we must have $N_i=N_0$ and $a_{N_i}^{(i)}=a_{N_0}^{(0)}$. Namely, Claim~\ref{cliam-4.1} holds.

\begin{clai}\label{cliam-4.2}
The situation $N\ge 3$ is impossible.
\end{clai}

Let $N\ge 3$. Taking $j = 1$ in \eqref{add-eqa}, we have
\begin{eqnarray}
\nonumber  Y_{i,1}(t) &\!\!\!=\!\!\!& \frac{1}{2}\big((t-i\alpha+2)Y_{0,1}(t) - (t-i\alpha)Y_{0,1}(t-i) - (t+\alpha)Y_{1,1}(t+1) \\
\label{deformation-eqa-j=1}   &\!\!\!\!\!\!& + i\alpha Y_{1,1}(t)+ (t-(i-1)\alpha)Y_{1,1}(t-i+1)\big).
\end{eqnarray}
Note that if we substitute \eqref{deformation-eqa-j=1} into \eqref{add-eqa}, then we obtain an equation on functions $Y_{0,1}$ and $Y_{1,1}$, and
both sides can be viewed as functions on $i,j$ and $t$.
Let ${\bf{L}}_2(i,j,t)$ and ${\bf{R}}_2(i,j,t)$ denote respectively the left- and right-hand sides of \eqref{add-eqa}.
By Claim~\ref{cliam-4.1}, a direct computation shows that there always exist $K_{N-2}(i,j), K_{0}(i,j)\in\C[i,j]$ such that
\begin{eqnarray*}
\mathrm{Coeff}_{{\bf{L}}_2(i,j,t) - {\bf{R}}_2(i,j,t)} t^{N-2} &\!\!\!=\!\!\!& \frac{1}{2}ijK_{N-2}(i,j), \\
\mathrm{Coeff}_{{\bf{L}}_2(i,j,t) - {\bf{R}}_2(i,j,t)} t^{0} &\!\!\!=\!\!\!& \frac{1}{2}ij\alpha K_{0}(i,j).
\end{eqnarray*}
By \eqref{add-eqa}, we have $K_{N-2}(i,j)=0$ if $ij\ne0$, and $K_{0}(i,j)=0$ if $ij\alpha\ne0$.
Through a more detailed analysis, we observe that $\deg\,K_{N-2}(i,j)\le 1$ and $\deg\,K_{0}(i,j)\le N-1$.
Since one can first fix $i\ne 0$ (resp., $j\ne 0$) and then let $j\rightarrow\infty$ (resp., $i\rightarrow\infty$), the above two observations imply that
(see Example~\ref{exam-N=3} for concrete equations in the case of $N=3$)

{\bf Case 1:} $\alpha \neq 0$.
\begin{equation*}
\left\{
\begin{aligned}
\mathrm{Coeff}_{K_{N-2}(i,j)}i &= 0,\\
\mathrm{Coeff}_{K_{N-2}(i,j)}j &= 0,\\
\mathrm{Coeff}_{K_{0}(i,j)}i^{N-1} &= 0,\\
\mathrm{Coeff}_{K_{0}(i,j)}j^{N-1} &= 0;
\end{aligned}
\right.
\end{equation*}

{\bf Case 2:} $\alpha = 0$.
\begin{equation*}
\left\{
\begin{aligned}
\left.\mathrm{Coeff}_{K_{N-2}(i,j)}i\right|_{\alpha = 0} &= 0,\\
\left.\mathrm{Coeff}_{K_{N-2}(i,j)}j\right|_{\alpha = 0} &= 0.
\end{aligned}
\right.
\end{equation*}
In both cases, one can derive that $a_{N} = 0$, which contradicts to our previous assumption \eqref{Yi1t}.
Hence, Claim~\ref{cliam-4.2} holds.

\begin{clai}\label{cliam-4.3}
The situation $N=2$ is impossible.
\end{clai}

Let $N=2$. Recall Claim~\ref{cliam-4.1}.
In this case, we have
\begin{equation*}
Y_{i,1}(t) = a_{2}t^2+ a_{1}^{(i)}t +a_{0}^{(i)},
\end{equation*}
which is essentially derived from relations \eqref{start-equations}.
Similarly, start from relations
\begin{equation}\label{start-equations-2}
[L_{j,2},L_{i-j,0}] = (i-2j)L_{i,2} -2 L_{i,3}, \quad [L_{-j,0},L_{j,2}]= 2jL_{0,2} + 2L_{0,3},
\end{equation}
one can derive that
\begin{equation*}
Y_{i,2}(t) = b_{2}t^2+ b_{1}^{(i)}t +b_{0}^{(i)}.
\end{equation*}
By comparing the coefficients of $t^2$ on both sides of \eqref{L01,Li0} with $i = 1$, we can derive that $a_{2} = 0$,
a contradiction. Hence, Claim~\ref{cliam-4.3} holds.

\begin{clai}\label{cliam-4.4}
If $N\le 1$, then ${Y_{i,m}(t)}=\beta_m(t-i\alpha) + \gamma_m$, where $\beta_m, \gamma_m\in\C$ and $\beta_0=1, \gamma_0=0$.
\end{clai}

Let $N\le 1$. Recall Claim~\ref{cliam-4.1}.
In this case, we have
\begin{equation*}
Y_{i,1}(t) = a_{1}t +a_{0}^{(i)}.
\end{equation*}
By \eqref{add-eqa-j=i}, we can derive that $a_{0}^{(i)}=a_{0}^{(0)}-i\alpha a_1$. Redenote $a_{1}$ by $\beta_1$, and $a_{0}^{(0)}$ by $\gamma_1$.
The above formula becomes
\begin{equation*}
Y_{i,1}(t) = \beta_1(t - i\alpha)+\gamma_1.
\end{equation*}
Generally, start from relations (cf.~\eqref{start-equations} for the case $m=1$, and \eqref{start-equations-2} for the case $m=2$)
\begin{equation*}\label{start-equations-m}
[L_{j,m},L_{i-j,0}] = (i-2j)L_{i,m} -m L_{i,m+1}, \quad [L_{-j,0},L_{j,m}]= 2jL_{0,m} + m L_{0,m+1},
\end{equation*}
one can derive that ${Y_{i,m}(t)}=\beta_m(t-i\alpha) + \gamma_m$ for some $\beta_m, \gamma_m\in\C$.
Finally, from \eqref{w-1-action-i0} we see that $\beta_0=1$ and $\gamma_0=0$.
This completes the proof of Claim~\ref{cliam-4.4}.

\begin{clai}\label{cliam-4.5}
We have $\beta_m = \beta^m$ and $\gamma_m = -m\alpha\beta^{m+1}$, where $\beta\in\C$.
\end{clai}

First, let us determine $\beta_m$. By Claim~\ref{cliam-4.4}, applying the relation
$$
[L_{0,m},L_{1,n}]=L_{1,m+n}+(n-m)L_{1,m+n+1}
$$
on $1$, one can obtain an equation on $\beta_m$ and $\gamma_m$.
By comparing the coefficients of $t$ on both sides of this equation, we obtain
$$
\beta_{m+n}+(n-m)\beta_{m+n+1} = \beta_m\beta_n + n\beta_m\beta_{n+1} - m\beta_n\beta_{m+1}.
$$
Recall that $\beta_0=1$ (by Claim~\ref{cliam-4.4}).
By Lemma~\ref{Sequence-problem}, we have $\beta_m = \beta^m$, where $\beta = \beta_1$.

Next, we determine $\gamma_m$. By Claim~\ref{cliam-4.4}, applying the relation $(i\ne 0)$
$$
[L_{0,m},L_{i,0}]=i L_{i,m}-m L_{i,m+1}
$$
on $1$, one can easily obtain $\gamma_m = -m\alpha\beta^{m+1}$.
This completes the proof of Claim~\ref{cliam-4.5}.

Finally, from Claims~\ref{cliam-4.1}--\ref{cliam-4.5} we see that this lemma holds.
\QED

\begin{exam}\label{exam-N=3}
\rm
For the convenience of the reader we explicitly write the systems of equations in Cases 1 and 2 in Claim~\ref{cliam-4.2} under the assumption $N=3$.
If $\alpha\ne0$, the system of equations is
\begin{equation*}
\left\{
\begin{aligned}
(1+2\alpha)\left(2(a_2^{(0)}-a_2^{(1)})-3(1+\alpha)a_3^{(0)}\right) &= 0,\\
(1+2\alpha)(a_2^{(0)}-a_2^{(1)})-(1+3\alpha+8\alpha^2)a_3^{(0)} &= 0,\\
2(a_2^{(0)}-a_2^{(1)})-3(1+\alpha)a_3^{(0)} &= 0,\\
(1-\alpha)(a_2^{(0)}-a_2^{(1)})-a_3^{(0)} &= 0;
\end{aligned}
\right.
\end{equation*}
and if $\alpha=0$, the system of equations is
\begin{equation*}
\left\{
\begin{aligned}
2(a_2^{(0)}-a_2^{(1)})-3a_3^{(0)} &= 0,\\
(a_2^{(0)}-a_2^{(1)})-a_3^{(0)} &= 0.
\end{aligned}
\right.
\end{equation*}
It is straightforward to check that both cases yield $a_3=a_3^{(0)}=0$.
\end{exam}

\subsection{\large{Simplicity and isomorphism classification of $\mathcal{W}(-1)$-modules}}

Similar to Theorem~\ref{w1-simplicity}, we have the following result on the simplicity of $\Omega_{\mathcal{W}(-1)}(\lambda,\alpha,\beta)$.

\begin{theo}\label{w-1-simplicity}\parskip-3pt
\begin{itemize}\parskip-3pt
\item[{\rm(1)}] $\Omega_{\mathcal{W}(-1)}(\lambda,\alpha,\beta)$ is simple if and only if $\alpha\neq0$.
\item[{\rm(2)}] $\Omega_{\mathcal{W}(-1)}(\lambda,0,\beta)$ has a unique proper submodule $t\Omega_{\mathcal{W}(-1)}(\lambda,0,\beta)\cong\Omega_{\mathcal{W}(-1)}(\lambda,1,\beta)$, and the quotient $\Omega_{\mathcal{W}(-1)}(\lambda,0,\beta)/t\Omega_{\mathcal{W}(-1)}(\lambda,0,\beta)$ is a one-dimensional trivial ${\mathcal{W}(-1)}$-module.
\end{itemize}
\end{theo}

\ni{\it Proof.}\ \
(1) Let $M=\Omega_{\mathcal{W}(-1)}(\lambda,\alpha,\beta)$. If $\alpha\neq0$, by viewing $M$ as a ${\rm Vir}$-module, from Lemma~\ref{lemm-13}(1), we see that $M$ is simple.
If $\alpha=0$, one can easily see that $t M$ is a submodule of $M$.

(2) From definition, one can easily see that $t\Omega_{\mathcal{W}(-1)}(\lambda,0,\beta)$ is the unique proper submodule of $\Omega_{\mathcal{W}(-1)}(\lambda,0,\beta)$.
Next, we prove $t\Omega_{\mathcal{W}(-1)}(\lambda,0,\beta)\cong\Omega_{\mathcal{W}(-1)}(\lambda,1,\beta)$ by comparing the actions of $\mathcal{W}(-1)$ on these two modules.

First, we consider the action of $\mathcal{W}(-1)$ on $t\Omega_{\mathcal{W}(-1)}(\lambda,0,\beta)$.
For $k\ge 0$, by \eqref{action-of-w-1}, we have
\begin{equation}\label{w-1-on-t0}
L_{i,m}\cdot (t\cdot t^k) = t\left(\sum_{s=0}^{k+1}(-1)^ss!\binom{m+s-1}{s}\binom{k+1}{s}\lambda^i\beta^{m+s}(t-i)^{k-s+1}\right).
\end{equation}
Next, we consider the action of $\mathcal{W}(-1)$ on $\Omega_{\mathcal{W}(-1)}(\lambda,1,\beta)$. Let $k\ge 0$.
By \eqref{action-of-w-1}, we have
\begin{eqnarray}
\nonumber L_{i,m}\cdot t^k &\!\!\!=\!\!\!& \sum_{s=0}^{k}(-1)^ss!\binom{m+s-1}{s}\binom{k}{s}\lambda^i\beta^{m+s}(t-i-(m+s)\beta)(t-i)^{k-s} \\
\nonumber &\!\!\!=\!\!\!& \sum_{s=0}^{k}(-1)^ss!\binom{m+s-1}{s}\binom{k}{s}\lambda^i\beta^{m+s}(t-i)^{k-s+1} \\
\nonumber &\!\!\!\!\!\!&- \sum_{s=0}^{k}(-1)^ss!\binom{m+s-1}{s}\binom{k}{s}\lambda^i\beta^{m+s+1}(m+s)(t-i)^{k-s} \\
\nonumber &\!\!\!=\!\!\!& \lambda^i\beta^m(t-i)^{k+1} + \sum_{s=1}^{k}(-1)^ss!\binom{m+s-1}{s}\binom{k}{s}\lambda^i\beta^{m+s}(t-i)^{k-s+1} \\
\nonumber &\!\!\!\!\!\!&- \sum_{s=0}^{k-1}(-1)^ss!\binom{m+s-1}{s}\binom{k}{s}\lambda^i\beta^{m+s+1}(m+s)(t-i)^{k-s} \\
\nonumber &\!\!\!\!\!\!&- (-1)^kk!\binom{m+k-1}{k}\lambda^i\beta^{m+k+1}(m+k) \\
\nonumber &\!\!\!=\!\!\!& \lambda^i\beta^m(t-i)^{k+1} + \sum_{s=1}^{k}(-1)^ss!\binom{m+s-1}{s}\binom{k}{s}\lambda^i\beta^{m+s}(t-i)^{k-s+1} \\
\nonumber &\!\!\!\!\!\!&+ \sum_{s=1}^{k}(-1)^s(s-1)!\binom{m+s-2}{s-1}\binom{k}{s-1}\lambda^i\beta^{m+s}(m+s-1)(t-i)^{k-s+1} \\
\nonumber &\!\!\!\!\!\!&+ (-1)^{k+1}(k+1)!\binom{m+k}{k+1}\lambda^i\beta^{m+k+1} \\
\nonumber &\!\!\!=\!\!\!& \lambda^i\beta^m(t-i)^{k+1} + \sum_{s=1}^{k}(-1)^ss!\binom{m+s-1}{s}\binom{k}{s}\lambda^i\beta^{m+s}(t-i)^{k-s+1} \\
\nonumber &\!\!\!\!\!\!&+ \sum_{s=1}^{k}(-1)^ss!\binom{m+s-1}{s}\binom{k}{s-1}\lambda^i\beta^{m+s}(t-i)^{k-s+1} \\
\nonumber &\!\!\!\!\!\!&+ (-1)^{k+1}(k+1)!\binom{m+k}{k+1}\lambda^i\beta^{m+k+1} \\
\label{w-1-on-1} &\!\!\!=\!\!\!& \sum_{s=0}^{k+1}(-1)^ss!\binom{m+s-1}{s}\binom{k+1}{s}\lambda^i\beta^{m+s}(t-i)^{k-s+1}.
\end{eqnarray}
Comparing \eqref{w-1-on-t0} with \eqref{w-1-on-1}, we see that $t\Omega_{\mathcal{W}(-1)}(\lambda,0,\beta)\cong\Omega_{\mathcal{W}(-1)}(\lambda,1,\beta)$.

At last, it is clear that the quotient $\Omega_{\mathcal{W}(-1)}(\lambda,0,\beta)/t\Omega_{\mathcal{W}(-1)}(\lambda,0,\beta)$ is a one-dimensional trivial ${\mathcal{W}(-1)}$-module.
This completes the proof.
\QED
\vskip15pt

Similar to Theorem~\ref{w1-iso-classification}, one can also prove the following result on the isomorphism classification of $\Omega_{\mathcal{W}(-1)}(\lambda,\alpha,\beta)$;
the details are omitted.

\begin{theo}\label{thm-42}
$\Omega_{\mathcal{W}(-1)}(\lambda_1,\alpha_1,\beta_1)\cong\Omega_{\mathcal{W}(-1)}(\lambda_2,\alpha_2,\beta_2)$
if and only if $\lambda_1=\lambda_2$, $\alpha_1=\alpha_2$ and $\beta_1=\beta_2$.
\end{theo}

\vskip10pt

\small \ni{\bf Acknowledgement}
We are very grateful to the referee for helpful comments that have improved the manuscript.
We also thank Professor Y. Su for helpful discussions.
This work was supported by the Fundamental Research Funds for the Central Universities (No. 2019QNA34).

\end{document}